\documentclass[a4paper,review]{cas-sc}

\usepackage[section]{placeins} 
\usepackage{amssymb}
\usepackage{amsthm}
\usepackage{amsmath}
\usepackage{upgreek}
\usepackage{pdflscape}
\usepackage{listings}
\usepackage{multirow}
\usepackage[percent]{overpic}
\usepackage{multicol}
\usepackage{color}
\usepackage{stmaryrd}
\usepackage{xfrac}
\usepackage[capitalise]{cleveref}
\usepackage{booktabs}
\usepackage[section]{placeins} 
\usepackage[section]{algorithm}
\usepackage{calc}
\usepackage[titletoc]{appendix}
\usepackage{siunitx}
\usepackage[sort&compress,square,numbers]{natbib}

\usepackage{graphicx}
\usepackage{graphics}
\usepackage{wrapfig}
\usepackage{float}
\usepackage{subfig}
\usepackage[percent]{overpic}
\usepackage{varwidth}
\usepackage{tikz}
\usepackage{pgfplots}
\usepackage{adjustbox}
\usetikzlibrary{arrows,matrix,positioning,fit}
\usetikzlibrary{shapes,positioning}
\usetikzlibrary{backgrounds}
\usepackage{tikz-layers}
\usepgfplotslibrary{fillbetween}
\usetikzlibrary{intersections}
\pgfplotsset{compat=1.14}
\usepgfplotslibrary{colorbrewer}
\usepgfplotslibrary{patchplots}
\usepgfplotslibrary[colorbrewer]
\usetikzlibrary{pgfplots.colorbrewer}
\usetikzlibrary[pgfplots.colorbrewer]
\usepgfplotslibrary{units}
\usetikzlibrary{spy}
\usepackage{pgfplotstable}
\usepackage{arrayjobx}
\graphicspath{ {./figs/} }
\usetikzlibrary{external}

\newlength\myheight
\newlength\mydepth
\settototalheight\myheight{Xygp}
\newcommand*\inlinegraphics[1]{%
  \settototalheight\myheight{Xygp}%
  \settodepth\mydepth{Xygp}%
  \raisebox{-\mydepth}{\includegraphics[height=\myheight]{#1}}%
}
\newcommand\orcid[1]{\href{https://orcid.org/#1}{\inlinegraphics{orcid_16x16.png}}}

\makeatletter
\def\BState{\State\hskip-\ALG@thistlm}
\makeatother

\newdefinition{definition}{Definition}[section]

\newcommand\pxvar[2]{\partial_{#2} #1}

\newcommand{\gu}{g(\mathbf{u})}

\newcommand{\gbar}{g(\overline{\mathbf{u}})}
\newcommand{\ux}{\mathbf{u}(\mathbf{x})}
\newcommand{\ubar}{\overline{\mathbf{u}}}

\newcommand{\ghat}{g(\widehat{\mathbf{u}}(\mathbf{x}))}

\begin{document}

\title[mode=title]{A note on higher-order and nonlinear limiting approaches for continuously bounds-preserving discontinuous Galerkin methods}
\shorttitle{Higher-order and nonlinear limiting approaches for continuously bounds-preserving discontinuous Galerkin methods}
\shortauthors{T. Dzanic}

\author[1]{T. Dzanic}[orcid=0000-0003-3791-1134]
\cormark[1]
\cortext[cor1]{Corresponding author}
\ead{dzanic1@llnl.gov}

\address[1]{Center for Applied Scientific Computing, Lawrence Livermore National Laboratory, Livermore, CA 94551, USA}

\begin{abstract}
In (Dzanic, \textit{J. Comp. Phys.}, 508:113010, 2024), a limiting approach for high-order discontinuous Galerkin schemes was introduced which allowed for imposing constraints on the solution \emph{continuously} (i.e., everywhere within the element). While exact for linear constraint functionals, this approach only imposed a sufficient (but not the minimum necessary) amount of limiting for nonlinear constraint functionals. This short note shows how this limiting approach can be extended to allow exactness for general nonlinear quasiconcave constraint functionals through a nonlinear limiting procedure, reducing unnecessary numerical dissipation. Some examples are shown for nonlinear pressure and entropy constraints in the compressible gas dynamics equations, where both analytic and iterative approaches are used. 
\end{abstract}

\begin{keywords}
Discontinuous Galerkin \sep
High-order\sep
Hyperbolic conservation laws\sep
Bounds-preserving\sep
Positivity-preserving\sep
Limiting
\end{keywords}



\maketitle

\section{Introduction}
\label{sec:intro}
The use of high-order discontinuous Galerkin (DG) schemes for simulating transport-dominated physics often requires additional numerical stabilization to ensure the discrete solution is \emph{structure-preserving} (i.e., is well-behaved and satisfies physical constraints). One such way of achieving this desirable behavior is applying some form of \emph{a posteriori} limiting on the solution, where the high-order DG approximation is blended with a secondary, more robust (but generally less accurate) approximation to yield a bounds-preserving (i.e., constraint-satisfying) solution, a few examples of which are presented in the following references~\citep{Zhang2010,Kuzmin2016,Anderson2017,Pazner2021,Hajduk2021,Dzanic2022,RuedaRamrez2022,Zhang2023,Lin2023,Peyvan2023,Ching2024}. This form of limiting has been widely used for simulating complex physical systems, particularly ones of predominantly hyperbolic nature, where numerical solutions must abide by known physical constraints (e.g., positivity of density and pressure in gas dynamics). However, a potential drawback of these limiting methods is that they are typically performed at discrete nodal locations for the DG solution, which only ensures that the given constraints are satisfied at those points. While adequate for many problems, the use of discretely bounds-preserving limiting may not be robust enough for applications where the DG solution must be evaluated at arbitrary locations (e.g., coupled meshes/solvers, remapping in arbitrary Lagrangian-Eulerian methods, adaptive mesh refinement, etc.). In such scenarios, the limited solution may still violate the constraints at these arbitrary points, causing the failure of the numerical scheme. 

A potential remedy for this problem was introduced in \citet{Dzanic2024} through a limiting approach that ensures the limited solution is \emph{continuously} bounds-preserving (i.e., across the entire solution polynomial). The approach relied on a novel functional for computing the limiting factor for the ``squeeze'' limiter of \citet{Zhang2010}, such that if the minimum of this limiting functional within an element was found (via a spatial optimization algorithm), the limited solution was guaranteed to be continuously bounds-preserving for any arbitrary quasiconcave constraint functional. This approach differed from other continuously bounds-preserving methods which rely on bounded basis functions (e.g., Bernstein polynomials)~\citep{Anderson2017, Glaubitz2019,Hajduk2021}, which may suffer from numerical inefficiencies and excessive dissipation (a more in-depth discussion on these issues is presented in \citet{Dzanic2024}, Section 1). The underlying mathematical basis of this approach was a linearization of the constraint functionals being enforced on the solution. For linear constraints (e.g., maximum principle on scalar solutions, positivity of density in gas dynamics, etc.), the computed limiting factor was ``exact'' (i.e., it was the minimum necessary amount of limiting). However, for nonlinear constraint functionals, the linearization could only ensure that the computed limiting factor was \emph{sufficient}, such that a smaller limiting factor (i.e., less numerical dissipation) could still yield a solution which satisfied the constraints continuously. Therefore, the limiting approach was suboptimal for nonlinear constraints commonly encountered in hyperbolic systems of equations as it was not applying only the minimum amount of limiting necessary. 

The purpose of this short note is to show how this limiting functional can be further modified such that the computed limiting factor is exact for arbitrary nonlinear quasiconcave functionals, reducing unnecessary numerical dissipation. While this proposed modification can be applied to general hyperbolic conservation laws, we focus specifically on the compressible Euler equations which admit higher-order and nonlinear constraints for pressure and entropy. Some examples are shown for the new approach, where both analytic and iterative approaches can be applied. 

\section{Preliminaries}
\label{sec:preliminaries}
In this section, some preliminaries are presented for the proposed modification, which briefly summarize the limiting approach introduced in \citet{Dzanic2024}. 

\subsection{Governing equations and constraints}
The exemplar hyperbolic conservation law used for this work is the compressible Euler equations for gas dynamics in $d$ dimensions, written in conservation form as
\begin{equation}\label{eq:gen_hype}
        \pxvar{\mathbf{u}(\mathbf{x}, t)}{t} + \boldsymbol{\nabla}\cdot\mathbf{F}(\mathbf{u}) = 0,
\end{equation}
where 
\begin{equation}\label{eq:euler}
    \mathbf{u} = \begin{bmatrix}
            \rho \\ \mathbf{m} \\ E
        \end{bmatrix} \quad  \mathrm{and} \quad \mathbf{F}(\mathbf{u}) = \begin{bmatrix}
            \mathbf{m}^T\\
            \mathbf{m}\otimes\mathbf{v} + P\mathbf{I}\\
        (E+P)\mathbf{v}^T
    \end{bmatrix}.
\end{equation}
Here, $\rho$ is the density, $\mathbf{m}$ is the momentum, and $E$ is the total energy. The symbol $\mathbf{I}$ denotes the identity matrix in $\mathbb R^{d \times d}$, $\mathbf{v} = \mathbf{m}/\rho$ denotes the velocity, and $P$ denotes the pressure, computed (with the assumption of a calorically perfect gas) as 
\begin{equation}
    P = (\gamma - 1)\rho e = (\gamma - 1)\left (E - \frac{1}{2} \mathbf{m}{\cdot}\mathbf{m}/ \rho \right),
\end{equation}
where $\gamma = 1.4$ is the specific heat ratio and $e$ is the specific internal energy. 

The solution of the Euler equations is endowed with a convex invariant set corresponding to the positivity of density ($\rho \geq 0$), positivity of pressure/internal energy ($P, \rho e \geq 0$), and a minimum principle on the specific physical entropy ($\sigma \geq \sigma_{\min}$), where $\sigma = P\rho^{-\gamma}$~\citep{Frid2001}. This can be represented by the positivity of a set of constraint functionals as 
\begin{equation}
    g_1(\mathbf{u}) = \rho - \rho_{\min}, \quad g_2(\mathbf{u}) = P - P_{\min}, \quad \text{and} \quad g_3(\mathbf{u}) = \sigma - \sigma_{\min},
\end{equation}
where small tolerances ($\rho_{\min}$, $P_{\min}$) are included for numerical stability purposes. This set of constraints is particularly interesting as they each have a unique mathematical character, consisting of a linear density constraint, a quadratic pressure/internal energy constraint, and a nonlinear entropy constraint. In many limiting approaches, the entropy constraint is neglected and only positivity of density/pressure is enforced (e.g., positivity-preserving limiters).

\subsection{Discontinuous Galerkin methods}
The underlying numerical method for this approach is the discontinuous Galerkin scheme (see, for example, \citet{Hesthaven2008DG}), where the solution $\mathbf{u}(\mathbf{x})$ within each element of the mesh is represented by a set of $n_s$ polynomial basis functions as 
\begin{equation}
    \mathbf{u}(\mathbf{x}) = \sum_{i=1}^{n_s} \mathbf{u}_{i}\phi_i(\mathbf{x})\subset V_h,
\end{equation}
where $\phi_i(\mathbf{x})$ are the basis functions, $\mathbf{u}_i$ are their associated coefficients, and $V_h$ is the piece-wise polynomial space spanned by the basis functions.
We consider a semi-discrete weak formulation of a hyperbolic conservation law in the form of 
\begin{equation}\label{eq:semi-disc}
    \sum_{k=1}^{N} \left \{ \int_{\Omega_k} \partial_t \mathbf{u}{\cdot} \mathbf{w}\ \mathrm{d}V  + \int_{\partial \Omega_k} \hat{\mathbf{F}}(\mathbf{u}^-, \mathbf{u}^+, \mathbf{n}) {\cdot} \mathbf{w}\ \mathrm{d}S  -  \int_{\Omega_k} \mathbf{F}(\mathbf{u}){\cdot} \nabla \mathbf{w}\ \mathrm{d}V\right \} = 0,
\end{equation}
where $\mathbf{w}(\mathbf{x}) \subset V_h$ is a test function and $\hat{\mathbf{F}}(\mathbf{u}^-, \mathbf{u}^+, \mathbf{n})$ is a numerical interface flux. Furthermore, we define the element-wise mean for an arbitrary element $\Omega_k$ as 
\begin{equation}
    \overline{\mathbf{u}}_k = \frac{\int_{\Omega_k}\mathbf{u} (\mathbf{x})\ \mathrm{d}\mathbf{x}}{\int_{\Omega_k} \mathrm{d}\mathbf{x}},
\end{equation}
which, under some relatively minor assumptions on the numerical scheme, generally preserves convex invariants of hyperbolic systems, such that it may be used as secondary bounds-preserving approximation for a limiting approach (see \citet{Zhang2010}, \citet{Zhang2011b}, and derived works).

\subsection{Continuously bounds-preserving limiting}
The proposed limiting approach in \citet{Dzanic2024} relies on the \emph{a posteriori} ``squeeze'' limiter of \citet{Zhang2010}, which linearly contracts the high-order DG solution towards the element-wise mean based on a limiting factor $\alpha$ as  
\begin{equation}\label{eq:limiter}
    \widehat{\mathbf{u}}(\mathbf{x}) = (1 - \alpha) \ux + \alpha \ubar = \ux + \alpha \left (\ubar - \ux \right).
\end{equation}
We drop the subscript $k$ for brevity as this limiter applies to every element in the domain. The goal of the proposed approach is to find a limiting factor such that 
\begin{equation}\label{eq:ghat}
    \ghat \geq 0 \ \forall \ \mathbf{x},
\end{equation}
for any arbitrary quasiconcave constraint functional $g(\mathbf{u})$. This approach differed from standard \emph{a posteriori} limiting methods for DG schemes in that it enforced constraints continuously instead of at discrete nodal locations (here, $\forall \ \mathbf{x}$ is used to refer to all locations \emph{within the element}). It was shown in \citet{Dzanic2024} that if one introduces the modified limiting functional $h(\mathbf{u})$ as 
\begin{equation}\label{eq:h}
    h(\mathbf{u}) = 
    \begin{cases}
    h^+(\mathbf{u}), \quad \text{if } \gu \geq 0, \\
    h^-(\mathbf{u}), \quad \text{else},
    \end{cases}
\end{equation}
where 
\begin{equation}\label{eq:hpm}
    h^+(\mathbf{u}) = \frac{\gu}{\gbar} \quad \text{and} \quad 
    h^-(\mathbf{u}) = \frac{\gu}{\gbar - \gu},
\end{equation}
then setting the limiting factor $\alpha$ as 
\begin{equation}\label{eq:alpha}
    \alpha = \max\left[0, -\underset{\mathbf{x}}{\min} \ h(\mathbf{u})\right],
\end{equation}
guarantees a limiting scheme that is continuously bounds-preserving.

In this approach, one can define ``exactness'' in the limiting as finding the minimum necessary value of $\alpha$ such that \cref{eq:ghat} is satisfied, which, for high-order solutions that initially violate the constraints, can be expressed as $\underset{\mathbf{x}}{\min}\ \ghat = 0$. However, for nonlinear quasiconcave constraint functionals, the linearization used in the definition of $h^-(\mathbf{u})$ can only ensure that the limiting is sufficient, i.e., $\underset{\mathbf{x}}{\min}\ \ghat \geq 0$.

\section{Proposed modifications}\label{sec:methodology}
The primary motivation of this work is to introduce a modification to the definition of $h^-(\mathbf{u})$ such that exactness can also be ensured for nonlinear functionals, reducing unnecessary numerical dissipation in the limiting procedure. Exactness for the limiting factor can be represented as finding the minimum necessary value of $\alpha$ such that constraints are satisfied, i.e.,
\begin{equation*}
    \underset{\alpha \geq 0}{\arg \min}\ g\left((1 - \alpha) \ux + \alpha \ubar \right) \geq 0\ \forall \ \mathbf{x},
\end{equation*}
which is equivalent to finding the spatial maximum of the necessary limiting factor at every point in the element, i.e.,
\begin{equation*}
    \underset{\mathbf{x}}{\max}\ \left[ 0, \alpha^*(\mathbf{x})\right ], \quad \text{where} \quad  g\left((1 - \alpha^*(\mathbf{x})) \mathbf{u}(\mathbf{x}) + \alpha^*(\mathbf{x}) \ubar \right) = 0.
\end{equation*}
To achieve this, we propose replacing the linearized formulation of $h^-(\mathbf{u})$ in \cref{eq:h}, denoted by $h^-_L(\mathbf{u})$, with a better (nonlinear) approximation of the necessary limiting factor, denoted by $h^-_{NL}(\mathbf{u})$, which mimics the above condition as 
\begin{equation}\label{eq:hnl}
    h_{NL}^-(\mathbf{u}(\mathbf{x})) = -\alpha^*(\mathbf{x}) :  \ghat = g\left((1 - \alpha^*(\mathbf{x})) \mathbf{u}(\mathbf{x}) + \alpha^*(\mathbf{x}) \ubar \right) = 0
\end{equation}
It can then be trivially shown following the proof in \citet{Dzanic2024} that setting the limiting factor as \cref{eq:alpha} using this modified formulation ensures $\underset{\mathbf{x}}{\min}\ \ghat = 0$ if the high-order solution is bounds-violating and $\ghat = \gu$ if the high-order solution is bounds-preserving, achieving exactness for arbitrary nonlinear constraint functionals. Furthermore, it also ensures $C^0$ continuity in $h(\mathbf{u})$ as $h_{NL}^-(\mathbf{u}) = h^+(\mathbf{u})$ when $g(\mathbf{u}) = 0$.

However, this formulation requires solving an intersection problem for the zero contour of the constraint functional. For linear constraints, it can be seen that this reduces to the linearized formulation, i.e., $h^-_{NL}(\mathbf{u}) = h^-_{L}(\mathbf{u})$. For quadratic constraints such as pressure, this higher-order limiting approach also admits an analytic solution, where $\alpha^*$ can be computed as the (positive) root of the quadratic equation
\begin{equation*}
   A \left(\alpha^*\right)^2 + B \alpha^* + C = 0,
\end{equation*}
with the coefficients taking on the values
\begin{equation*}
    A = \Delta \rho \Delta E - \frac{1}{2} \Delta \mathbf{m} \cdot \Delta \mathbf{m}, 
    \quad B = E \Delta \rho + \rho \Delta E - \mathbf{m} \cdot \Delta \mathbf{m} - \Delta \rho \frac{P_{\min}}{\gamma - 1}, 
    \quad C = \rho E - \frac{1}{2} \mathbf{m} \cdot \mathbf{m} - \rho \frac{P_{\min}}{\gamma - 1},
\end{equation*}
for positive pressure constraints, where $\Delta \mathbf{u} = \ubar - \mathbf{u}$. This quadratic formulation is similarly seen in the works of \citet{Kuzmin2016} and \citet{Zhang2011b}. Similar closed-form expressions can be used for higher-order constraints (e.g., cubic) which admit analytic solutions to the intersection problem in \cref{eq:hnl}. 

If one also wants to ensure exactness for nonlinear constraints such as entropy which do not have an analytic solution to the intersection problem, it can be straightforwardly solved using simple iterative root-bracketing approaches such as the bisection or Illinois method, which is similar to discrete limiting approaches with entropy-based constraints~\citep{Guermond2019}. The quasiconcavity of the constraint functionals ensures that this root-bracketing problem is well-behaved, i.e., the is solution is bounded by $0 \leq \alpha^* \leq 1$ and is unique. Furthermore, one can apply a stricter upper bound for the initial bracket by using the linearized formulation in \cref{eq:hpm},  further reducing the computational cost.

\section{Implementation}
\label{sec:implementation}

The modified functional was implemented within the same numerical framework as in \citet{Dzanic2024} and closely follows the original work, using identical parameters and optimization approaches. This approach is summarized as follows, but for more information, the reader is referred to Section 3 of \citet{Dzanic2024}. The solution was first converted into a modal basis in monomial form for fast evaluation at arbitrary locations. The initial guess for constraint functional for the optimization process was evaluated at the solution nodes (e.g., Gauss--Lobatto nodes) as well as the volume/surface quadrature nodes (e.g., Gauss--Legendre nodes), the minimum of which was taken as the initial point for optimization. We note that the distinction between solution nodes and quadrature nodes is not necessarily present in collocated flux reconstruction schemes with closed solution nodes, but we evaluate at these nodes regardless to mimic standard nodal DG approaches. Optimization was performed using two iterations of the Newton--Raphson method with a fallback on adaptive gradient descent with backtracking line search in regions where Newton--Raphson is ill-suited. The Jacobian and Hessian of the constraint functionals were computed numerically, and a bound on the true minimum of the constraint functional was extrapolated at the end of the optimization steps. For the gas dynamics equations, the constraints were enforced sequentially, first on density, then pressure, then entropy, to ensure that the constraint functionals were well-behaved. Minimum density and pressure values were set as $\rho_{\min} = P_{\min} = 10^{-11}$. A minimum tolerance was set as $\epsilon = 10^{-12}$, and the limiting factor was explicitly set as $\alpha = 1$ when $g(\ubar) < \epsilon$ to avoid numerically undefined behavior.

Some further implementation details are necessary for the proposed modification to the approach. For pressure constraints, it is necessary to compute $h_{NL}^-(\mathbf{u})$ as the root of the quadratic polynomial. It can be shown that in locations where $h(\mathbf{u}) = h_{NL}^-(\mathbf{u})$ (i.e., where $g(\mathbf{u}) < 0$)  and under the assumption that $g(\ubar) > 0$, the discriminant of the quadratic equation is positive and only one of the roots, given by 
\begin{equation*}
    \frac{-B + \sqrt{B^2 - 4AC}}{2A},
\end{equation*}
yields a limiting factor $\alpha^* \in [0,1]$. In practice, care must be taken to avoid numerical round-off errors in calculating this root, namely by enforcing a floor value of zero for the discriminant. The edge case of $A \sim 0$ is taken care of in a similar manner to the edge case of $g(\ubar) \sim 0$, where the limiting factor is explicitly set as $\alpha^* = 1$ when $|A| < \epsilon$. Therefore, for computational purposes, the limiting factor is computed as
\begin{equation}
    \alpha^* = 
    \begin{cases}
    1, \hspace{82pt} \text{if } |A| < \epsilon, \\
    \frac{-B + \sqrt{\max ( 0, B^2 - 4AC)}}{2A}, \quad \text{else}.
    \end{cases}
\end{equation}

For the entropy constraints, 5 iterations of the Illinois method were used to compute $h_{NL}^-(\mathbf{u})$ in locations where $g(\mathbf{u}) < 0$. The initial lower (bounds-violating) bracket was set as $\alpha_l = 0$ with $g_l = g(\alpha_l)$, and the upper (bounds-preserving) bracket was set as $\alpha_h = -h_{L}^-(\mathbf{u})$ with $g_h = g(\alpha_h)$ (recall that the linearized formulation always yields an upper bound on the exact limiting factor necessary). At each iteration, the next guess was computed as 
\begin{equation}
    \alpha_m = \frac{\alpha_l g_h - \alpha_h g_l}{g_h - g_l},
\end{equation}
and the brackets were updated as 
\begin{equation*}
    \alpha_l = \alpha_m, \quad g_l = g(\alpha_m), \quad \text{and} \quad g_h = \frac{1}{2}g(\alpha_h)
\end{equation*}
if $g(\alpha_m) < 0$ and 
\begin{equation*}
    \alpha_h = \alpha_m, \quad g_l = \frac{1}{2}g(\alpha_l), \quad \text{and} \quad g_h = g(\alpha_m)
\end{equation*}
if $g(\alpha_m) \geq 0$. At the end of the root-finding process, the solution was set to the upper bound, i.e., $\alpha^* = -\alpha_h$, to ensure that the bounds-preserving estimate of the limiting factor was chosen. 
\section{Results}\label{sec:results}

To highlight the improvements of the proposed modification, we consider first a pathological example of limiting a static discontinuity within one element, where a discontinuity is placed in the center of the element with the left/right states set as 
\begin{equation}
    \left[\rho, u, P \right]^T = \begin{cases}
        \left[1, 1, 2 P_{\min} \right]^T, &\mbox{if } x \leq 0.5, \\
        \left[3, 3, 1 \right]^T, &\mbox{else}.
    \end{cases} 
\end{equation}
A 9th-order polynomial DG approximation was initialized by interpolating the solution on the Gauss--Lobatto nodes spaced along the domain $\Omega = [0,1]$, and constraints were enforced for positivity of density/pressure as well as a minimum entropy principle. An arbitrary minimum entropy value was set as $\sigma_{\min} = 0.1$. These conditions yield various properties for the unlimited solution to showcase the behavior of the new approach, with a continuously bounds-preserving density field, a discretely (but not continuously) bounds-preserving pressure field, and a discretely bounds-violating entropy field. The minimum entropy value was chosen to demonstrate the nonlinear limiting approach for bounds-violating solution as well as to mimic a condition encountered in time-dependent simulations of hyperbolic conservation laws, where entropy bounds (which might be computed from the previous time step) are violated by the discrete solution at the next time step~\citep{Dzanic2022, Ching2024}.

    \begin{figure}[htbp!]
        \centering
        \subfloat[Pressure]{\adjustbox{width=0.4\linewidth,valign=b}{    \begin{tikzpicture}[spy using outlines={rectangle, height=3cm,width=2.3cm, magnification=3, connect spies}]
		\begin{axis}[name=plot1,
		    axis x line=left,
            axis y line=left,
		    xlabel={$x$},
		    xtick={0,0.2,0.4,0.6,0.8,1},
    		xmin=0,
    		xmax=1,
    		x tick label style={
        		/pgf/number format/.cd,
            	fixed,
            	fixed zerofill,
            	precision=1,
        	    /tikz/.cd},
    		ylabel={$P$},
    		ylabel style={rotate=-90},
		    ytick={-1, -0.5, 0, 0.5, 1, 1.5},
    		ymin=-1,
    		ymax=1.5,
    		y tick label style={
        		/pgf/number format/.cd,
            	fixed,
            	fixed zerofill,
            	precision=1,
        	    /tikz/.cd},
    		legend style={at={(0.97, 0.03)},anchor=south east ,font=\small, column sep=0.2cm},
    		legend cell align={left},
    		style={font=\normalsize}]

            \addplot[mark=none, color={blue!80}, style={thick}] coordinates {(0,0) (0.5, 0.0) (0.5, 1.0) (1.0, 1.0)};
    		\addlegendentry{Exact solution}    		

			\addplot[color={black}, style={dashed, thick}]
				table[x=x,y=P,col sep=comma,unbounded coords=jump]{./figs/data/example_linear_continuous.csv};
    		\addlegendentry{Unlimited solution}    		

			\addplot[color={red!90!black}, style={dashed}]
				table[x=x,y=Pl,col sep=comma,unbounded coords=jump]{./figs/data/example_linear_continuous.csv};
    		\addlegendentry{Linear limiting}
      
			\addplot[color={red!90!black}, style={}]
				table[x=x,y=Pl,col sep=comma,unbounded coords=jump]{./figs/data/example_nonlinear_continuous.csv};
    		\addlegendentry{Nonlinear limiting}

			\addplot[color=black, style={ultra thin}, only marks, mark=o, mark options={scale=0.6}]
				table[x=x,y=P,col sep=comma,unbounded coords=jump]{./figs/data/example_linear_nodes.csv};
      
            \addplot[mark=none, color={black!80}, style={dotted, thick}] coordinates {(0,0) (1, 0)} ;

		\end{axis}

	\end{tikzpicture}}}
        \subfloat[Entropy]{\adjustbox{width=0.4\linewidth,valign=b}{    \begin{tikzpicture}[spy using outlines={rectangle, height=3cm,width=2.3cm, magnification=3, connect spies}]
		\begin{axis}[name=plot1,
		    axis x line=left,
            axis y line=left,
		    xlabel={$x$},
		    xtick={0,0.2,0.4,0.6,0.8,1},
    		xmin=0,
    		xmax=1,
    		x tick label style={
        		/pgf/number format/.cd,
            	fixed,
            	fixed zerofill,
            	precision=1,
        	    /tikz/.cd},
    		ylabel={$\sigma$},
    		ylabel style={rotate=-90},
		    ytick={-0.5, 0, 0.5},
    		ymin=-0.5,
    		ymax=0.5,
    		y tick label style={
        		/pgf/number format/.cd,
            	fixed,
            	fixed zerofill,
            	precision=1,
        	    /tikz/.cd},
    		legend style={at={(0.97, 0.03)},anchor=south east ,font=\small, column sep=0.2cm},
    		legend cell align={left},
    		style={font=\normalsize}]
    		
            \addplot[mark=none, color={blue!80}, style={thick}] coordinates {(0,0) (0.5, 0.0) (0.5, 0.21479) (1.0, 0.21479)};

			\addplot[color={black}, style={thick}]
				table[x=x,y=s,col sep=comma,unbounded coords=jump]{./figs/data/example_linear_continuous.csv};

			\addplot[color={red!90!black}, style={dashed}]
				table[x=x,y=sl,col sep=comma,unbounded coords=jump]{./figs/data/example_linear_continuous.csv};
      
			\addplot[color={red!90!black}, style={}]
				table[x=x,y=sl,col sep=comma,unbounded coords=jump]{./figs/data/example_nonlinear_continuous.csv};

			\addplot[color=black, style={ultra thin}, only marks, mark=o, mark options={scale=0.6}]
				table[x=x,y=s,col sep=comma,unbounded coords=jump]{./figs/data/example_linear_nodes.csv};
      
            \addplot[mark=none, color={black!80}, style={dotted, thick}] coordinates {(0,0.1) (1, 0.1)};

		\end{axis}

	\end{tikzpicture}}}
        \caption{\label{fig:static_example} Comparison of the linear and nonlinear limiting 
        approaches for pressure (left) and entropy (right) constraints for a static discontinuity problem computed with a $\mathbb P_9$ approximation. Pressure limiting computed using analytic formulation, entropy limiting computed using iterative formulation. Dotted line represents zero contour for the constraint functional.}
    \end{figure}
As the density field was already continuously bounds-preserving, the limiter was first applied to the pressure field. A comparison of the linearized formulation and the proposed nonlinear formulation is shown in \cref{fig:static_example}. It can be seen that the proposed nonlinear formulation drastically reduces the numerical dissipation in the scheme, with the limited solution significantly closer to the unlimited solution in comparison to the linearized formulation. Furthermore, it can be seen that the new limiting approach is exact in the sense that the minimum pressure of the limited solution was on the order of $P_{\min}$. Similar results could be seen for the entropy constraints, also shown in \cref{fig:static_example}. Due to the stronger nonlinearity of the constraint functional, the linearized approach showed even more unnecessary dissipation than with the pressure constraint, whereas the proposed nonlinear approach was exact and showed a significant reduction in the numerical dissipation. 

A more quantitative evaluation of the proposed modification was performed with a near-vacuum isentropic Euler vortex which possesses an analytic solution for comparison. The problem consists of a smooth vortex on a periodic domain $\Omega = [-10, 10]^2$, with the initial solution set as 
    \begin{equation*}
        \begin{bmatrix}
        \rho \\
        u \\
        v \\
        P
        \end{bmatrix}
        = \begin{bmatrix}
            P^\frac{1}{\gamma} \\
            \frac{S}{2 \pi R} (y-y_0)\phi(\mathbf{x}) \\
            1 - \frac{S}{2 \pi R} (x-x_0)\phi(\mathbf{x}) \\
            \frac{1}{\gamma M^2} \left(1 - \frac{S^2 M^2 (\gamma-1)} {8 \pi^2}\phi(\mathbf{x})^2\right)^\frac{\gamma}{\gamma-1}
        \end{bmatrix}, \quad \mathrm{where} \quad \phi(\mathbf{x}) = \exp{\left(\frac{1-\|\mathbf{x}-\mathbf{x}_0\|_2^2}{2R^2}\right)}.
    \end{equation*}
The parameters are set to $R = 1.5$, $M = 0.4$, and $S = 28.11711$, which yield a near-vacuum state at the vortex peak with a minimum density of $\rho = 8{\cdot}10^{-9}$ and minimum pressure of $P = 2{\cdot}10^{-11} = 2P_{\min}$. One flow-through of the domain ($t = 20$) was computed using both a $\mathbb P_4$ and $\mathbb P_5$ approximation using meshes of varying resolution with only positivity-preserving constraints, after which the $L^{\infty}$ norm of the pressure error (computed at solution nodes) was compared from the solutions obtained using the original linear limiting approach and the proposed nonlinear limiting approach. The comparisons for both the $\mathbb P_4$ and $\mathbb P_5$ approximations are shown in \cref{tab:icv4} and \cref{tab:icv5}, respectively. It can be seen that the proposed modification significantly reduces the overall error in the pressure field, with both approximation orders showing decreases between $20-70\%$ across the varying levels of mesh resolution. These results highlight the reduced numerical dissipation that stems from the modified constraint functional proposed in this note. Furthermore, the imposition of exact pressure constraints (as opposed to a linearized approximation) can have further benefits in terms of the computational cost as the reduced pressure (i.e., reduced speed of sound) can increase the maximum admissible explicit time step for near-vacuum problems. A simple computational cost comparison (at a fixed time step) was performed for the near-vacuum isentropic Euler vortex problem using a $\mathbb P_4$ approximation with $N = 30^2$ elements on one NVIDIA V100 GPU in terms of the absolute wall clock time. As the added effort of analytically computing $h_{NL}^-(\mathbf{u})$ was essentially negligible compared to computing $h_{L}^-(\mathbf{u})$, the computational cost of the approach with the proposed modification was effectively identical to the linearized method, both with a computational cost increase of approximately $18\%$ over a discrete limiting approach. When $h_{NL}^-(\mathbf{u})$ was instead computed numerically (which is not necessary for the pressure constraints used for this example but is simply presented for comparison), the overall computational cost increase in comparison to a discrete limiting method was then $24.7\%$. We remark here that the relatively marginal difference between the analytic and numerical nonlinear approaches can largely be attributed to the efficacy of GPU computing for local compute-intensive operations.  

    \begin{figure}[htbp!] 
        \centering
        \begin{tabular}{lcccccc}
        \toprule
        Method & $N = 20^2$ & $N = 30^2$ & $N = 40^2$ & $N = 50^2$ & $N = 60^2$ & $N = 70^2$ \\ 
        \midrule
Linear limiting & $2.86\times 10^{-1}$ & $2.93\times 10^{-2}$ & $1.61\times 10^{-2}$ & $7.17\times 10^{-3}$ & $3.40\times 10^{-3}$ & $2.34\times 10^{-3}$ \\
Nonlinear limiting & $1.55\times 10^{-1}$ & $2.36\times 10^{-2}$ & $9.88\times 10^{-3}$ & $4.33\times 10^{-3}$ & $1.61\times 10^{-3}$ & $6.27\times 10^{-4}$ \\
\midrule
    Error reduction & \textcolor{green!70!black}{-45.9\%} & \textcolor{green!70!black}{-19.4\%} & \textcolor{green!70!black}{-38.5\%} & \textcolor{green!70!black}{-39.6\%} & \textcolor{green!70!black}{-52.8\%} & \textcolor{green!70!black}{-73.2\%} \\
        \bottomrule
        \end{tabular}
        \captionof{table}{\label{tab:icv4} Comparison of the $L^{\infty}$ norm of the pressure error using linear and nonlinear limiting for the near-vacuum isentropic Euler vortex after one flow-through of the domain ($t = 20$) computed with a $\mathbb P_4$ approximation. }
    \end{figure}
    
    \begin{figure}[htbp!] 
        \centering
        \begin{tabular}{lcccccc}
        \toprule
        Method & $N = 20^2$ & $N = 30^2$ & $N = 40^2$ & $N = 50^2$ & $N = 60^2$ & $N = 70^2$ \\ 
        \midrule
        Linear limiting & $1.15\times 10^{-1}$ & $2.02\times 10^{-2}$ & $3.58\times 10^{-3}$ & $3.89\times 10^{-3}$ & $1.97\times 10^{-3}$ & $4.18\times 10^{-4}$ \\
        Nonlinear limiting & $4.97\times 10^{-2}$ & $1.66\times 10^{-2}$ & $2.78\times 10^{-3}$ & $1.24\times 10^{-3}$ & $5.36\times 10^{-4}$ & $2.43\times 10^{-4}$ \\
        \midrule
        Error reduction & \textcolor{green!70!black}{-56.5\%} & \textcolor{green!70!black}{-21.4\%} & \textcolor{green!70!black}{-27.1\%} & \textcolor{green!70!black}{-70.1\%} & \textcolor{green!70!black}{-72.0\%} & \textcolor{green!70!black}{-37.7\%} \\
        \bottomrule
        \end{tabular}
        \captionof{table}{\label{tab:icv5} Comparison of the $L^{\infty}$ norm of the pressure error using linear and nonlinear limiting for the near-vacuum isentropic Euler vortex after one flow-through of the domain ($t = 20$) computed with a $\mathbb P_5$ approximation. }
    \end{figure}
\section{Concluding remarks}\label{sec:conclusion}
We proposed an improvement to the continuously bounds-preserving limiting approach presented in \citet{Dzanic2024}, which allows for the approach to achieve the exact amount of limiting necessary for arbitrary nonlinear quasiconcave constraint functionals as opposed to just a sufficient amount. The modification relies on replacing the linearization in limiting functional with an intersection/root-finding problem, which may be computed analytically for some constraints and numerically for others. Some examples were shown for the compressible Euler equations, showing the reduced numerical dissipation and increased accuracy of the proposed approach.

\section*{Acknowledgements}
\label{sec:ack}

This work was performed under the auspices of the U.S. Department of Energy by Lawrence Livermore National Laboratory under contract DE--AC52--07NA27344 and the LLNL-LDRD Program under Project tracking No.\ 24--ERD--050. Release number LLNL--JRNL--863292.

\bibliographystyle{unsrtnat}
\bibliography{reference}

\begin{thebibliography}{17}
\providecommand{\natexlab}[1]{#1}
\providecommand{\url}[1]{\texttt{#1}}
\expandafter\ifx\csname urlstyle\endcsname\relax
  \providecommand{\doi}[1]{doi: #1}\else
  \providecommand{\doi}{doi: \begingroup \urlstyle{rm}\Url}\fi

\bibitem[Zhang and Shu(2010)]{Zhang2010}
Xiangxiong Zhang and Chi-Wang Shu.
\newblock On maximum-principle-satisfying high order schemes for scalar conservation laws.
\newblock \emph{Journal of Computational Physics}, 229\penalty0 (9):\penalty0 3091--3120, May 2010.
\newblock \doi{10.1016/j.jcp.2009.12.030}.

\bibitem[Kuzmin and Lohmann(2016)]{Kuzmin2016}
Dmitri Kuzmin and Christoph Lohmann.
\newblock \emph{Synchronized slope limiting in discontinuous Galerkin methods for the equations of gas dynamics}.
\newblock Technische Universit{\"a}t Dortmund, Fakult{\"a}t f{\"u}r Mathematik, 2016.

\bibitem[Anderson et~al.(2017)Anderson, Dobrev, Kolev, Kuzmin, de~Luna, Rieben, and Tomov]{Anderson2017}
R.~Anderson, V.~Dobrev, Tz. Kolev, D.~Kuzmin, M.~Quezada de~Luna, R.~Rieben, and V.~Tomov.
\newblock High-order local maximum principle preserving ({MPP}) discontinuous {G}alerkin finite element method for the transport equation.
\newblock \emph{Journal of Computational Physics}, 334:\penalty0 102--124, April 2017.
\newblock \doi{10.1016/j.jcp.2016.12.031}.

\bibitem[Pazner(2021)]{Pazner2021}
Will Pazner.
\newblock Sparse invariant domain preserving discontinuous {G}alerkin methods with subcell convex limiting.
\newblock \emph{Computer Methods in Applied Mechanics and Engineering}, 382:\penalty0 113876, August 2021.
\newblock \doi{10.1016/j.cma.2021.113876}.

\bibitem[Hajduk(2021)]{Hajduk2021}
Hennes Hajduk.
\newblock Monolithic convex limiting in discontinuous {G}alerkin discretizations of hyperbolic conservation laws.
\newblock \emph{Computers \& Mathematics with Applications}, 87:\penalty0 120–138, April 2021.
\newblock \doi{10.1016/j.camwa.2021.02.012}.

\bibitem[Dzanic and Witherden(2022)]{Dzanic2022}
T.~Dzanic and F.D. Witherden.
\newblock Positivity-preserving entropy-based adaptive filtering for discontinuous spectral element methods.
\newblock \emph{Journal of Computational Physics}, 468:\penalty0 111501, November 2022.
\newblock \doi{10.1016/j.jcp.2022.111501}.

\bibitem[Rueda-Ram{\'{\i}}rez et~al.(2022)Rueda-Ram{\'{\i}}rez, Pazner, and Gassner]{RuedaRamrez2022}
Andr{\'{e}}s~M. Rueda-Ram{\'{\i}}rez, Will Pazner, and Gregor~J. Gassner.
\newblock Subcell limiting strategies for discontinuous {G}alerkin spectral element methods.
\newblock \emph{Computers \& Fluids}, 247:\penalty0 105627, October 2022.
\newblock \doi{10.1016/j.compfluid.2022.105627}.

\bibitem[Zhang and Cheng(2023)]{Zhang2023}
Fan Zhang and Jian Cheng.
\newblock Analysis on physical-constraint-preserving high-order discontinuous {G}alerkin method for solving {K}apila's five-equation model.
\newblock \emph{Journal of Computational Physics}, 492:\penalty0 112417, November 2023.
\newblock \doi{10.1016/j.jcp.2023.112417}.

\bibitem[Lin et~al.(2023)Lin, Chan, and Tomas]{Lin2023}
Yimin Lin, Jesse Chan, and Ignacio Tomas.
\newblock A positivity preserving strategy for entropy stable discontinuous {G}alerkin discretizations of the compressible {E}uler and {N}avier-{S}tokes equations.
\newblock \emph{Journal of Computational Physics}, 475:\penalty0 111850, February 2023.
\newblock \doi{10.1016/j.jcp.2022.111850}.

\bibitem[Peyvan et~al.(2023)Peyvan, Shukla, Chan, and Karniadakis]{Peyvan2023}
Ahmad Peyvan, Khemraj Shukla, Jesse Chan, and George Karniadakis.
\newblock High-order methods for hypersonic flows with strong shocks and real chemistry.
\newblock \emph{Journal of Computational Physics}, 490:\penalty0 112310, October 2023.
\newblock \doi{10.1016/j.jcp.2023.112310}.

\bibitem[Ching et~al.(2024)Ching, Johnson, and Kercher]{Ching2024}
Eric~J. Ching, Ryan~F. Johnson, and Andrew~D. Kercher.
\newblock Positivity-preserving and entropy-bounded discontinuous {G}alerkin method for the chemically reacting, compressible {E}uler equations. {P}art i: The one-dimensional case.
\newblock \emph{Journal of Computational Physics}, 505:\penalty0 112881, May 2024.
\newblock \doi{10.1016/j.jcp.2024.112881}.

\bibitem[Dzanic(2024)]{Dzanic2024}
Tarik Dzanic.
\newblock Continuously bounds-preserving discontinuous {G}alerkin methods for hyperbolic conservation laws.
\newblock \emph{Journal of Computational Physics}, 508:\penalty0 113010, July 2024.
\newblock \doi{10.1016/j.jcp.2024.113010}.

\bibitem[Glaubitz(2019)]{Glaubitz2019}
Jan Glaubitz.
\newblock Shock capturing by {B}ernstein polynomials for scalar conservation laws.
\newblock \emph{Applied Mathematics and Computation}, 363:\penalty0 124593, December 2019.
\newblock \doi{10.1016/j.amc.2019.124593}.

\bibitem[Frid(2001)]{Frid2001}
Hermano Frid.
\newblock Maps of convex sets and invariant regions for finite-difference systems of conservation laws.
\newblock \emph{Archive for Rational Mechanics and Analysis}, 160\penalty0 (3):\penalty0 245--269, November 2001.
\newblock \doi{10.1007/s002050100166}.

\bibitem[Hesthaven and Warburton(2008)]{Hesthaven2008DG}
Jan~S. Hesthaven and Tim Warburton.
\newblock \emph{Nodal Discontinuous {G}alerkin Methods}.
\newblock Springer New York, 2008.
\newblock \doi{10.1007/978-0-387-72067-8}.

\bibitem[Zhang and Shu(2011)]{Zhang2011b}
Xiangxiong Zhang and Chi-Wang Shu.
\newblock On positivity-preserving high order discontinuous {G}alerkin schemes for compressible {E}uler equations on rectangular meshes.
\newblock \emph{Journal of Computational Physics}, 229\penalty0 (23):\penalty0 8917--8934, November 2011.
\newblock \doi{10.1016/j.jcp.2010.08.016}.

\bibitem[Guermond et~al.(2019)Guermond, Popov, and Tomas]{Guermond2019}
Jean-Luc Guermond, Bojan Popov, and Ignacio Tomas.
\newblock Invariant domain preserving discretization-independent schemes and convex limiting for hyperbolic systems.
\newblock \emph{Computer Methods in Applied Mechanics and Engineering}, 347:\penalty0 143--175, April 2019.
\newblock \doi{10.1016/j.cma.2018.11.036}.

\end{thebibliography}



\end{document}